\documentclass[11pt]{article}

\usepackage[T1]{fontenc}
\usepackage[english]{babel}
\usepackage{amsmath}
\usepackage{amsthm}
\usepackage{amssymb}
\usepackage{amsfonts}
\usepackage{latexsym}
\usepackage{euscript}
\usepackage{multicol}
\usepackage[dvips]{graphicx}
\newtheorem{theo}{Theorem}[section]
\newtheorem{Def}{Definition}[section]
\newtheorem{Lem}{Lemma}[section]
\newtheorem{Rem}{Remark}[section]
\newtheorem{Prop}{Proposition}[section]
\newtheorem{Cor}{Corollary}[section]

\title{\Large Dimensional properties of the harmonic measure for a random walk on  a hyperbolic group}
\author{Vincent Le Prince}

\begin{document}
\maketitle

\noindent \begin{bf}Abstract:\end{bf}
This paper deals with random walks on isometry groups of Gromov hyperbolic spaces, and more precisely with the dimension of the harmonic measure $\nu$ associated with such a random walk. We first establish a link of the form $\dim \nu\leq h/l$ between the dimension of the harmonic measure, the asymptotic entropy $h$ of the random walk and its rate of escape $l$. Then we use this inequality to show that the dimension of this measure can be made arbitrarily small and deduce a result on the type of the harmonic measure.

\vspace{0.5cm} 

\noindent \begin{bf}Keywords:\end{bf} ergodic theory, random walk, hyperbolic group, harmonic measure, entropy.

\noindent \begin{bf}MSC:\end{bf} 60J15, 20F67, 28D20, 28A78.

\section*{Introduction}

%\ \quad 
Let $(X,d)$ be a hyperbolic space and $G$ a non-elementary subgroup of its isometry group acting properly discontinuously. Given a  probability measure $\mu$ on $G$, we define the associated random walk $(x_n)$ by $x_n=h_1\cdots h_n$, $h_i$ being $\mu$-distributed independent random variables. Under some hypotheses on $\mu$, we can describe the asymptotic behaviour of this random walk: if the support of $\mu$ generates $G$ then almost all trajectories converge to an element $x_{\infty}$ of the hyperbolic boundary $\partial X$ of $X$, whose distribution is called the harmonic measure and denoted by $\nu$; if moreover  $\mu$ has a finite first moment then $(\partial X,\nu)$ coincides with the Poisson boundary of the random walk (\cite{K00}, see Section \ref{Preliminaries and notations} for the definitions and details).
 
On the boundary of a tree we can define a metric by :
$$ dist(\xi_1 , \xi_2) = e^{-(\xi_1 | \xi_2)}\ ;$$
$(\xi_1 | \xi_2)$ being the Gromov product (here the length of the common part) of $\xi_1$ and $\xi_2$. In the general case, this formula doesn't define a metric, but there exists a family of metrics $(d_a)_{1<a<a_0}$ on  $\partial X$ having the same properties. 
We shall consider in this article the relations between this metric structure and the harmonic measure. More precisely we are  interested in the Hausdorff dimension of this measure and hence  in its type. Our main result will be the fact that the harmonic measure is not necessary absolutely continuous with respect to the Hausdorff measure on the boundary but can be singular. This question has been considered for random walks on $SL(2,\mathbb{R})$ related to continued fractions (\cite{Ch}, \cite{KPW}) and in some others contexts (e.g. \cite{PUZ}).

Note that in some situations like for the nearest neighbor random walks on finitely generated free groups (see \cite{DM}, \cite{LM}, \cite{L}) or on finite free products of finite groups (\cite{MM}), this harmonic measure can be explicitely computed, but in general this is not the case. Whereas in the above cited articles  combinatorial methods, based on a description of the hyperbolic boundary as a set of infinite words,  are used, we shall use  geometrical methods. 

\vspace{5pt}

Let us now describe more precisely the contents of this paper.\\
The following relation between the dimension of $\nu$, the rate of escape $l$ of the random walk and its asymptotic entropy $h$ (see Section \ref{Preliminaries and notations} for the definitions) was established in the case of trees in \cite{K98} and \cite{L} :
\begin{equation}\label{dim nu =h/l (eq)}
\dim \nu = \frac{h}{l}\ .
\end{equation}
Section \ref{Lien entre les quantité asymptotiques et la dimension de Hausdorff de la mesure harmonique} is devoted to an extension of this result to the case of the subgroup of isometries of a hyperbolic group  (see  Prop. \ref{majoration de la dimension ponctuelle par h/l} for a more precise statement): 
\begin{equation}\label{dim nu leq h/l (eq)}\dim \nu \leq \frac{1}{\log a}\ \frac{h}{l} \ ,\end{equation}
  $a$ being the parameter used for the choice of the metric on the hyperbolic boundary.\\
Then we want to show that the dimension of this measure can be arbitrarily small. For this purpose we construct in Section \ref{Construction} a sequence $(\mu_k)$ of probability measures such that, denoting the asymptotic entropy and the rate of escape associated with $\mu_k$  respectively by $h(G,\mu_k)$ and $l(G,\mu_k)$, we have 
$$\lim_{k\rightarrow +\infty}\frac{h(G,\mu_k)}{l(G,\mu_k)}= 0\ .$$
To show this property we will need a formula to estimate the rate of escape; this formula is stated in Section \ref{Une formule intégrale pour la vitesse de fuite}.
Combined with the previous upper bound on  the dimension, this gives us the following result.

\vspace{5pt}

\noindent\begin{bf}Theorem  \ref{gros théorème}\end{bf}
{\em Let $G$ be a subgroup of the group of isometries of a hyperbolic space $(X,d)$ acting properly discontinuously and which is not elementary. For every $\epsilon>0$ there exists on $G$ a symmetric probability measure $\mu$ with a finite first moment, whose support generates $G$, and which has the following property: the pointwise dimension of the harmonic measure $\nu$ associated with $(G,\mu)$ is $\nu$-almost surely smaller than $\epsilon$.}

\vspace{5pt}

\noindent In the case where $X$ is the Cayley graph of the hyperbolic group $G$ equipped with the word metric (w.r.t.\ a certain system of generators $S$), the dimension of the boundary $\partial G$ is equal to the growth $v(G,S)$ of $G$ (w.r.t.\ the same $S$). Moreover the Hausdorff measure is then finite and non zero. Theorem \ref{gros théorème} implies that the harmonic measure can be singular with respect to this Hausdorff measure (Corollary \ref{corollaire sur les groupes (Cor)}). 

\vspace{5pt}

Two other natural questions on the dimension are the following.\\
The first, which in some sense is the opposite to the one we consider, is the question of knowing whether the harmonic measure can be of maximal dimension. In the case where $X$ is the Cayley graph of $G$, the support of the harmonic measure is $\partial X$.  The dimension of $\partial X$ is, denoting by $v$ the growth of $G$ and working with the distance $d_a$ on $\partial X$,  equal to $v_a(G)=v/\log a$ (\cite{Coo2}, see Section \ref{Preliminaries and notations}). We have a fundamental relation between $h$, $l$ an $v$ (see \cite{V00}):
\begin{equation}\label{h leq l v (eq)}
h\leq l v.
\end{equation} 
In view of the inequality (\ref{dim nu leq h/l (eq)}), the question of the maximality of the dimension of the harmonic measure is then related to the question raised in \cite{V00} of knowing whether the quotient $h/l$ can be maximal, i.e. equal to $v$. 

\noindent The second  is the question of the positivity of the dimension of the harmonic measure. Under the hypothesis stated in the first paragraph of this introduction, the Poisson boundary is not trivial. A well-known fact of the theory of Poisson boundaries is that this non-triviality is equivalent to the asymptotic entropy $h$ being strictly positive (\cite{KV}). So in the case of a tree the formula $\dim\nu =h/l$ implies that the dimension of the harmonic measure is strictly positive. But with only the inequality (\ref{dim nu leq h/l (eq)}) we can't conclude. In \cite{KL} a similar result (i.e. positivity of the dimension of the harmonic measure) is proven in the case of the Brownian motion on a Riemannian manifold with pinched negative sectional curvature.

\section{Preliminaries and notations}\label{Preliminaries and notations}

\subsection{ Hyperbolic spaces}\label{espaces hyperboliques}
We will need a rather complete description of some aspects of the
geometry of hyperbolic spaces (\cite{G}), which we recall there. Our main references are  \cite{CDP} and \cite{GH}.

\subsubsection*{Definition}
Let $(X,d)$ be a proper geodesic metric space.  Let $o$ be a point in $X$; we
define the Gromov product (w.r.t.  $o$) on $X^2$ :
$$ (x|y)_o=\frac{1}{2}[d(o,x)+d(o,y)-d(x,y)]\ .$$
We call $(X,d)$ a $\delta$-hyperbolic space if this product
satisfies, for all $x$, $y$, $z$, $o$ in $X$,
\begin{equation}\label{hyperbolicité}
 (x|z)_o\geq \min\{(x|y)_o,(y|z)_o\}-\delta \ .
\end{equation}
We say that $X$ is  hyperbolic if there exists $\delta$ such that $X$ is $\delta$-hyperbolic.
An important class of hyperbolic spaces is the one of word hyperbolic groups. One can define on the Cayley graph  of a finitely generated group $G$ (w.r.t.\ a symmetric system of generators $S$)  a metric which is the word metric (w.r.t.\ $S$) in restriction to $G$ and which makes every edge isometric to the segment $[0,1]$. A finitely generated group $G$ is said to be hyperbolic when its Cayley graph equipped with this metric is a hypaerbolic space (which doesn't depend on the choice of $S$).

\noindent We fix for the sequel a $\delta$-hyperbolique space $(X,d)$ equipped with a point $o$ and we denote the Gromov product of $x$ and $y$ w.r.t. $o$ by $(x|y)$.

\subsubsection*{Hyperbolic boundary}
We recall that a geodesic (resp. geodesic ray, resp. geodesic segment)
in $X$ is an isometry from $\mathbb{R}$ (resp. $[0, \infty[$, resp.
$[a,b]$) to $X$, as well as the image of such an isometry. 

\noindent Two geodesic rays $\sigma_1$ and $\sigma_2$ in $X$ are said to be equivalent ($\sigma_1 \sim \sigma_2$)
if there exists a constant $D$ such that for all $t$,
$d(\sigma_1(t),\sigma_2(t))\leq D$. The hyperbolic boundary is defined as the quotient of the set of geodesic rays by this equivalence relation. 

\noindent One can  extend the Gromov product to $X\cup \partial X$.
% : for $x$, $y$
%in $\partial X$,
%\begin{equation}\label{Produit de Gromov sur le bord}
%(x|y)=\inf \{\liminf_n (x_n|y_n)\} \ ;
%\end{equation}
%the infimum being taken among the sequences $(x_n)$ and $(y_n)$
%which converge
% respectively to $x$ and $y$. 
The relation
(\ref{hyperbolicité}) remains true. One can then define a topology on $X\cup \partial X$, which makes $X\cup \partial X$
be a compactification of  $X$, by taking for each  point $\xi$ in $\partial X$ as
a base of neighborhood the sets
$$\big\{x\in X\cup \partial X : (\xi|x)>R\big\}\ ,\ (R>0).$$
A  sequence $(x_n)$ converges to an element $\xi$ in $\partial X$
if and only if
$$\lim_{n\rightarrow\infty}(\xi|x_n)=+\infty \ .$$
In particular, if $\sigma$ is a geodesic ray, $\sigma(t)$ converges, when $t$ goes to infinity, to the equivalence class of $\sigma$ in $\partial X$, which we denote by $\sigma(+\infty)$.
Moreover, if $(x_n)$ converges to $x$ and $(y_n)$ to $y$, we have (see \cite{GH})
\begin{equation}\label{}
(x|y) \leq \liminf_n (x_n|y_n)\leq (x|y) + 2\delta \ .
\end{equation}

\noindent The action of the isometry group on $X$ extends to a continuous action on the boundary.

\noindent When $X$ is the Cayley graph of a hyperbolic group $G$, we shall write $\partial G$ instead of $\partial X$.

\subsubsection*{Metric on the boundary}
On the boundary of a tree we can define a metric by :
$$ dist(\xi_1 , \xi_2) = e^{-(\xi_1 | \xi_2)}\ ;$$
In the general case, this formula doesn't define a metric, but there exists a family of metrics $(d_a)_{1<a<a_0}$ on  $X\cup\partial X$ (see \cite{CDP}, ch. 11). 
We fix such a $a$ for the sequel. We will need the following
property, which shows the analogy with the case of the trees
(see \cite{CDP}):
\begin{Prop}\label{propriété de la métrique sur le bord}
There exists a constant $\lambda >0$ such that
\begin{itemize}
\item $\forall x,y \in X\cup \partial X$, $d_a(x,y)\leq \lambda a^{-(x|y)}$ ;
\item $\forall x,y \in \partial X$, $d_a(x,y)\geq \lambda^{-1} a^{-(x|y)}$ .
\end{itemize}
\end{Prop}
\noindent Note that the topology induced by this metric  on $X\cup
\partial X$ coincides with the one previously defined.

\subsubsection*{Limit set of a subgroup of isometries}
Let $G$ a subgroup of isometries acting properly discontinuously
on $(X,d)$. A reference for this part is \cite{Coo1}.
\begin{Def}\label{ensemble limite}
The set of accumulation points in $\partial X$ of the $G$-orbit of
a point $x$ in $X$ doesn't depend on the point $x$. We call it the
limit set of the group $G$ denoted by $\Lambda_G$.
\end{Def}
\noindent In the case of a hyperbolic group  $G$ acting on its Cayley graph, $\Lambda_G=\partial G$. This limit set will be the support of the harmonic measure; we need some
of its properties.
\begin{Def}
The group $G$ is said to be elementary if $\Lambda_G$ consists of at
most two elements.
\end{Def}
\noindent If $G$ is non-elementary, $\Lambda_G$ is uncountable.

\begin{Def}\label{élément hyperbolique}
Let $g$ be an isometry in $X$. This element is said to be
hyperbolic if $(g^n o)_{n\in\mathbb{Z}}$ is a quasi-geodesic ;
which means that there exists two constants $\lambda$ and $c$ such
that for all $n$ and $m$ in $\mathbb{Z}$,
$$\lambda^{-1}|n-m|-c \leq d(g^n o,g^m o)\leq \lambda|n-m|+c \ .$$
\end{Def}
\begin{Prop}\label{quasi geodesique}
Let $\eta$ be a $(\lambda,c)$-quasi geodesic. Then $\eta(t)$ admits limits $\eta(\pm\infty)$ when $t$ goes to  $\pm\infty$. If $\sigma$ is a geodesic with endpoints $\sigma(\pm\infty)=\eta(\pm\infty)$, then the images of $\sigma$ and $\eta$ lie at a finite Hausdorff distance $K$ from each other, $K$ depending only on $\lambda$, $c$ and $\delta$.
\end{Prop}

\noindent In particular if an isometry $g$ is hyperbolic, then there exists $g^+$ and
$g^-$ in $\partial X$ such that $(g^n o)$ goes to $g^+$ and $g^-$
when $n$ go respectively to plus or minus infinity. In fact for all
$x$ in $X$, the sequences $(g^n x)$ converge to the same limits.
The points $g^+$ and $g^-$ are  fixed points for $g$, which are
said to be respectively attractive and repulsive.

\begin{Prop}(\cite{Coo1})
A non-elementary group G contains hyperbolic elements.
\end{Prop}
\noindent An important property of the limit set is its minimality (see
\cite{Coo1}).
\begin{Prop}(Gromov)
Assume $G$ to be non-elementary. Then every non-void $G$-invariant
compact set in $\partial X$ contains $\Lambda$. In other words,
$\Lambda$ is the only minimal set.
\end{Prop}
\noindent This has as a consequence :
\begin{Prop}\label{conséquence de la minimalité de l'ensemble limite}
Assume $G$ to be non-elementary. Let $U$ be an open set in
$\partial X$ which intersects $\Lambda$. Then we have :
$$\partial X = \bigcup_{g\in G}gU \ ;$$
and since $\partial X$ is compact : there exists a finite set of elements $g_1,\cdots
,g_r$ in $G$ such that
\begin{equation}\label{}
\partial X =  \bigcup_{1\leq i\leq r}g_i U \ .
\end{equation}
\end{Prop}

\subsection{Random walk}\label{Random walk}
We fix in this part a subgroup $G$ of the group of isometries of $(X,d)$ acting properly discontinuously and which is non-elementary.
\subsubsection*{Definition}
 Let $\mu$ be a probability measure on $G$. The random walk on $G$ determined by the measure $\mu$ is the Markov chain $\pmb{x}=(x_n)_{n\geq 0}$ with transition probabilities
$$p(x,y)=\mu(x^{-1}y)$$
and starting at $x_0=e$. The position $x_n$ at time $n\geq 1$ of the random walk is given by
$$x_n=h_1\cdots h_n\ ;$$  $(h_n)_{n\geq 1}$ being a sequence  of independent $\mu$-distributed random variables. We call the $h_n$ the increments and $(x_n)$ the trajectory of the random walk.
  We note $\mathbb{P}$ the distribution of $\pmb{x}$ in $G^{\mathbb{N}}$ and $T$ the shift in
$G^{\mathbb{N}}$. 

\noindent The Markov operator $P$ associated to the random walk is then:
$$Pf(x)= \sum_{g\in G}\mu(g)f(xg)\ .$$
We say that a function $f$ on $G$ is $\mu$-harmonic if it satisfies $Pf=f$.

\subsubsection*{Poisson boundary and harmonic measure}\label{Poisson boundary and harmonic measure}

The behaviour of the paths of the random walk is described by the following:
\begin{Prop}\label{prop : comportement et frontière de Poisson de la marche}(\cite{K00})
Assume that the support of the measure $\mu$ generates $G$. Then the random walk $(x_n o)$ associated to $\mu$ converges $\mathbb{P}$-almost surely to an element $x_{\infty}$ in $\partial X$. 
\end{Prop}
\noindent We denote by $bnd$ the map (defined on a set of $\mathbb{P}$-measure one) from $G^{\mathbb{N}}$ to $\partial X$ which associate $x_{\infty}$ to $\pmb{x}=(x_n)$ and we note $\nu$ the distribution of $x_{\infty}$, which we call the harmonic measure. Note that we have  $\nu=bnd(\mathbb{P})$ and that $bnd$ is G-equivariant.

\noindent Let us recall that if $G$, equipped with a measure $\mu$, acts on
a space $Y$ with a measure $m$, the convolution $\mu*m$ of $m$ by
$\mu$ is defined by : for all continuous bounded function $f$ ,
$$ \int_Y f(y)d(\mu*m)(y)=\int_{G\times Y}f(gy)d\mu(g)dm(y) \ .$$
We say that a measure is $\mu$-stationary if $\mu*\nu=\nu$.

%A topological space on which $G$ acts continuously is called a $G$-space. A $\mu$-boundary is a $G$-space equipped with a $\mu$-stationary probability measure $\lambda$ such that $(x_n \lambda)$ converges weakly to a $\delta$-measure for $\mathbb{P}$-almost all trajectory of the random walk $(x_n)$ associated with $(G,\mu)$. 

\noindent The measure $\nu$ is $\mu$-stationary. %doing from $(\partial X,\nu)$ a $\mu$-boundary. 
If we assume in addition that $\mu$ has a finite first moment,  that  is to say
$$\sum_{g\in G}\mu(g)d(o,go)<+\infty\ ,$$
then $(\partial X,\nu)$ is in fact the Poisson boundary associated to $(G,\mu)$(\cite{K00}); which means that every  $\mu$-harmonic bounded function $f$ on $G$ can be written :
$$f(g)=\int_{\partial X} F(g\xi) d\nu (\xi) \ ;$$
$F$ being a bounded function on $\partial X$ (for details on this notion, see \cite{F73} and \cite{KV}).
%maximal among the $\mu$-boundaries (\cite{K00}) : we call it the Poisson boundary and $\nu$ the harmonic measure. This name is due to the fact

\noindent A first result on the type of the harmonic measure is that under previous hypothesis it doesn't have any atom. Besides if the support of $\mu$ generates $G$ as a semi-group the support of the measure $\nu$ is the limit set $\Lambda_G$ of $G$. This is a consequence of the minimality of this limit set and of  stationnarity of $\nu$, and can be established using  Proposition
\ref{conséquence de la minimalité de l'ensemble limite} through the following lemma:
\begin{Lem}\label{la mesure harm charge les ouverts}
Let $U$ be an open set in $\partial X$ which intersects $\Lambda_G$. Then $\nu(U)>0$.
\end{Lem}
\begin{proof}
Let $\gamma_1,\cdots ,\gamma_r \in G$ be such that
$$\Lambda_G \subset \bigcup_{1\leq i \leq r}\gamma_i U \ ;$$
and choose $s$ such that $\gamma_1,\cdots ,\gamma_r \in supp(\mu^s)$. Since $\nu$ is $\mu$-stationnary we have
$$\nu(U)=\sum_{g\in supp(\mu^s)}\mu^s(g)\nu(gU) \ .$$
Hence if $\nu(U)=0$, we would have $\nu(\Lambda_G)=0$.
\end{proof}

\subsubsection*{Asymptotic quantities}\label{Asymptotic quantities}
Let $\mu$ be a probability measure on $G$ having a finite first moment.\\
We write $|g|=d(o,go)$, and for a probability measure $\lambda$,
$$L(\lambda):=\sum_{g\in G}\lambda(g) |g|  \ .$$

\noindent Denote by $\mu^n$ the $n^{\textrm{th}}$ convolution of $\mu$, which
is the distribution of the position at time $n$ of the random
walk. The sequence $(L(\mu^n))_n$ is subadditive, which allows us
to adopt the following definition    :
\begin{Def}
The limit of the  sequence $\big(L(\mu^n)/n\big)_n$
%$$l(G,\mu)=\lim_{n\rightarrow + \infty} \frac{1}{n}L(\mu_n)  \ $$
is called rate of escape of the random walk  $(x_n)_{n\geq 0}$ and denoted by $l(G,\mu)$.
\end{Def}
\noindent  Moreover,  $\mu$ having a finite first moment implies that the 
 entropy $$H(\mu):=\sum_g  -\mu(g)\log(\mu(g))$$  of $\mu$ is finite.
One can define in the same way the asymptotic entropy :
\begin{Def}
The limit of the  sequence $\big(H(\mu^n)/n\big)_n$
is called asymptotic entropy of the random walk  $(x_n)_{n\geq 0}$ and denoted by $h(G,\mu)$.
%The quantity
%$$h(G,\mu)=\lim_{n\rightarrow + \infty} \frac{1}{n}H(\mu^n)\ $$
%is called  asymptotic entropy of the random walk   $(x_n)_{n\geq 0}$
\end{Def}

\noindent Using Kingman's subbadditive ergodic theorem, one gets
(\cite{D})  a $\mathbb{P}$-almost sure convergences:
%(an analogousof Shannon-MacMilan-Breiman theorem) :
$$-\frac{1}{n} \log(\mu^n(x_n))\longrightarrow h(G,\mu) \quad\textrm{and}\quad
\frac{|x_n|}{n} \longrightarrow l(G,\mu)\ .$$
If there is no ambiguity, we shall denote by $h$ and $l$ these two
quantities. If in addition the support of $\mu$ generates $G$, both are positive (see \cite{K00}).

\subsection{Dimension}\label{Dimension}

Let $(X,d)$ be a complete metric space.  One defines the $\alpha$-Hausdorff measure of a set $Z\subset X$  as 
$$m_H(Z,\alpha)=\lim_{\epsilon\rightarrow 0} 
    \inf\Big\{\sum_{U\in \mathcal{G}_{\epsilon}}\big(diam (U)\big)^{\alpha}\Big\} \ ;$$ 
the infimum being taken over the covers $\mathcal{G}_{\epsilon}$ of $Z$ by open sets of diameter at most $\epsilon$.
The usual Hausdorff dimension of $Z$ is then 
\begin{equation}\label{définition de la dimension de H}
\dim_H(Z)=\inf\big\{\alpha : m_H(Z,\alpha)=0\big\}=\sup\big\{\alpha : m_H(Z,\alpha)=+\infty \big\} \ ;
\end{equation}
When $m_H(X,\dim_H(X))$ is finite and non-zero, the function $Z\mapsto m_H(Z,\dim_H X)$ is, after normalization, a probability measure on $X$, called the Hausdorff measure.

\subsubsection*{Dimensions of measures}
Let $\nu$ be a  probability measure on $(X,d)$. We define the Hausdorff dimension of $\nu$ as :
$$\dim_H \nu = \inf\{\dim_H Z : \nu(Z)=1\}\ .$$
\begin{Rem}\label{la dimension caractérise le type}
The dimension of a measure characterizes to some extent its type. Indeed if $\lambda$ is absolutely continuous w.r.t. $\nu$, $\{Z : \nu(Z)=1\}\subset\{ Z : \lambda(Z)=1\}$, so 
${\dim_H \lambda \leq \dim_H \nu}$.
\end{Rem}

\noindent To estimate the dimension of a set, and a fortiori of a measure, is usually not easy. We will have however a more direct way to estimate the dimension of a  measure, introducing the (lower and upper) pointwise dimensions at a point $x$ :
$$\underline{\dim}_P \nu (x)= \liminf_{r\rightarrow 0}\frac{\log \nu B(x,r)}{\log r}\ \ \  \textrm{and}\ \ \
                                     \overline{\dim}_P \nu (x)= \limsup_{r\rightarrow 0}\frac{\log \nu B(x,r)}{\log r} \ .$$
Moreover this notion allows a more intuitive vision of the dimension of a measure : it can be regarded as the rate of decrease of the measure of balls. 
To relate the Hausdorff and pointwise dimensions, we will need a condition on the space (see \cite{P}, appendix 1).
\begin{Def}
A metric space $(X,d)$ is said to have finite multiplicity if there  exist $K>0$ and $\epsilon_0 >0$ such that for all
$\epsilon \in ]0,\epsilon_0[$, there exists a cover of multiplicity $K$ (i.e. in which every point belongs at most to $K$ balls) of $X$ by balls of radius $\epsilon$. 
\end{Def}

\begin{Prop}\label{majoration de la dim de H par la dim ponctuelle}(see \cite{P}, appendix 1)
Let $(X,d)$  be a finite multiplicity space. Let $\nu$ be a  probability measure on $X$. If there exists a constant 
$d$ such that $\overline{\dim}_P \nu (x)\leq d $ $\nu$-almost surely, then ${\dim}_H \nu \leq d\ .$
\end{Prop}

\noindent If $\underline{\dim}_P \nu (x)\geq d $ $\nu$-almost surely, we also have
$\dim_H \nu \geq d$. A  probability measure with a constant $d$ such that $\underline{\dim}_P \nu (x) = \overline{\dim}_P \nu (x)= d$
$\nu$-almost surely is said to be exact dimensionnal. All dimensions are equal in this case.

\subsubsection*{Case of the boundary of a hyperbolic group}\label{Cas du bord d'un groupe hyperbolique}
We adopt the notations and hypotheses of the part \ref{Poisson boundary and harmonic measure}. As we want to consider the question of the maximality of the dimension of the harmonic measure, it is natural to ask about the dimension of the limit set  $\Lambda_G$, which supports this measure.
We define the critical exposant of base $a$ of $G$:
$$e_a(G)= \limsup_{R\rightarrow\infty}\frac{\log_a card\{g\in G : d(o,go)\leq R\}}{R}\ .$$
 
\begin{Prop}\label{dimension de l'ensemble limite d'un groupe q-c-c}(\cite{Coo2})
Under the hypothesis previously adopted and if in addition G is quasi-convex-cocompact (see \cite{Coo2} for a definition), then $\Lambda_G$ has Hausdorff dimension $e_a(G)$ (with respect to the metric $d_a$).
\end{Prop}
\noindent We don't want to explore this hypothesis of quasi-convex-cocompacity. We just note that in the case of a discrete group acting on the hyperbolic half-plane it coincides with the notion of convex-cocompacity; and we shall also use the fact that a hyperbolic group $G$ acting on its Cayley graph has this property (\cite{Coo2}).

\noindent In Section \ref{Lien entre les quantité asymptotiques et la dimension de Hausdorff de la mesure harmonique}, we bound from above the pointwise dimension of the harmonic measure. In order to obtain a result on the singularity of this measure with respect to the $e_a(G)$-Hausdorff measure, we need, in view of Proposition \ref{majoration de la dim de H par la dim ponctuelle} and Remark \ref{la dimension caractérise le type}, a finite multiplicity result on $(\partial X, d_a)$. Let us remark that when $X$ is a tree, this is obvious because balls of the same diameter form a partition.

\begin{theo}\label{multiplicité finie du bord d'un groupe hyperbolique}
We assume that $G$ is a hyperbolic group ($X$ its Cayley graph). Then its hyperbolic boundary $\partial G$, equipped with the metric $d_a$, is of finite multiplicity.
\end{theo}
\begin{proof}
We first note that if two geodesic rays joining $o$ with two points $\xi_1$ and $\xi_2$ of the boundary go through a same point $x$ such that $d(o,x)=n$, then  $(\xi_1|\xi_2)\geq n$ and so (see Proposition \ref{propriété de la métrique sur le bord})
$d_a (\xi_1,\xi_2)\leq \lambda a^{-n}$. Denote by $W_n$ the set of words $w$ of length $n$ through which goes a certain geodesic ray starting at $o$ ; and for each $w$ in $W_n$ let $\xi_w$ be the limit point of such a ray. In view of our first remark, the set of all open balls $B(\xi_w,\lambda a^{-(n-1)})$, $w$ in $W_n$, is a cover of the boundary. \\
%($n-1$ because we consider open balls).\\
 We are now going to show that these covers are of finite multiplicity (uniformly bounded in $n$).
Take $n>0$ and set $\epsilon = a^{-n}$. Let $\xi$ be a point in $\partial X$ and let $w$ be a  word of  length $n$ through which goes a ray $[o,\xi[$. Now let  $B(\xi_{w'},\epsilon)$ be a ball of our cover in which lies $\xi$. 
Since $d_a(\xi,\xi_{w'})<\epsilon$, $\lambda^{-1}a^{-(\xi|\xi_{w'})}<\lambda a^{-n}$ ; and so
$$(\xi|\xi_{w'})\geq n - 2\log_a \lambda \ .$$
But we have
$$(w|w')\geq \min\big\{(w|\xi),(\xi|\xi_{w'}),(\xi_{w'}|w')\big\} - 2 \delta\ ;$$
which yieds, as $(w|\xi)=(\xi_{w'}|w')=n$,
\begin{align*}
(w|w')&\geq \min\big\{n,(\xi|\xi_{w'})\big\} - 2 \delta\\
      & \geq n - 2\log_a \lambda  - 2 \delta \ .
\end{align*}
Since $(w|w')=n-\frac{1}{2}d(w,w')$, we deduce that
$$d(w,w')\leq  4\log_a \lambda  + 4 \delta \ .$$
But $G$ is finitely generated, so there is only a finite number of $w'$ of length $n$ which distance from $w$ is less than this constant, and this number doesn't depend on $n$.
\end{proof}
%\begin{Rem}\label{}
%This property extends to the  general case when $X$ is not a Cayley graph, but we need $G$ being quasi-convex-cocompact. Indeed if we assume $G$ having this property, then the $G$-orbit  of $o$ is quasi isometric to the Cayley graph of $G$ (see \cite{Coo2}).
%\end{Rem}

\section{Relation between asymptotic quantities and pointwise dimension of the harmonic measure}\label{Lien entre les quantité asymptotiques et la dimension de Hausdorff de la mesure harmonique}

In this section $(X,d)$ is a $\delta$-hyperbolic space, $G$ a non-elementary subgroup of the group of isometries on $X$ acting properly discontinuously.  We also fix  a measure $\mu$ on $G$ with a finite first moment, and such that its support generates $G$ as a semigroup. We shall use the notations of the previous section.

\noindent We have (see Section \ref{Asymptotic quantities}) $\mathbb{P}$-almost surely
$$\frac{(x_{n-1}|x_n)}{n}\longrightarrow l \quad \textrm{and} \quad-\frac{\log \mu^n(x_n)}{n}\longrightarrow h \ .$$
We define for $\epsilon>0$ and $N$  the set $\Omega_{\epsilon}^N$ of trajectories such that for $n\geq N$,
\begin{itemize}
\item{$(x_{n-1}|x_n)>(l-\epsilon)n$}
\item{$-\log\mu^n(x_n)<(h-\epsilon)n$.}
\end{itemize}
For $\eta>0$ there exists then an integer $N_{\epsilon,\eta}$ such that
$$\forall  N\geq N_{\epsilon,\eta}\ ,\ \mathbb{P}(\Omega_{\epsilon}^N)>1-\eta \ ;$$
we denote $\Omega_{\epsilon,\eta}=\Omega_{\epsilon}^{N_{\epsilon,\eta}}$.  Denote also by $C_{\pmb{x}}^n$ the set of trajectories whose $n^{\textrm{th}}$ position coincides with $x_n$. Our demonstration is based on the following lemma :
\begin{Lem}(\cite{K98})\label{Kaimanov. (Lem)}
There exist a set $\Lambda_{\epsilon,\eta}\subset\Omega_{\epsilon,\eta}$ with measure greater than $1-2\eta$, on which the quantity $$\frac{\mathbb{P}(C_{\pmb{x}}^n\cap\Omega_{\epsilon,\eta})}{\mu^n(x_n)}$$
admits a strictly positive limit when $n$ goes to infinity. In particular on this set we have
$$\limsup_n\frac{\log\mathbb{P}(C_{\pmb{x}}^n\cap\Omega_{\epsilon,\eta})}{n}
            =\limsup_n\frac{\log\mu^n(x_n)}{n} \ .$$
\end{Lem}
\noindent We write, for $\xi\in\partial X$ and $r>0$,
$$D(\xi,r)=\big\{\pmb{x}: x_{\infty}\in B(\xi,r)\big\}\ .$$
In the next lemma we prove that if a trajectory is at time $n$ in the same place as $x_n$, then its endpoint is not too far from $x_{\infty}$.
\begin{Lem}\label{not too far from x{infini} (Lem)}
Fix two strictly positive number $\eta$ and $\epsilon$. We have, for $\pmb{x}\in\Lambda_{\epsilon,\eta}$ and $n\geq N_{\epsilon,\eta}$,
$$C_{\pmb{x}}^n\cap\Omega_{\epsilon,\eta}\subset
    D\Bigg(x_{\infty},\frac{2\lambda a^{-(l-\epsilon)(n+1)}}{1-a^{-(l-\epsilon)}}\Bigg)\ ;$$
$\lambda$ being the constant introduced in Proposition \ref{propriété de la métrique sur le bord}.
\end{Lem}
\begin{proof}We fix $\pmb{x}\in\Lambda_{\epsilon,\eta}$.
Let $\pmb{x}'$ be an element of $C_{\pmb{x}}^n\cap\Omega_{\epsilon,\eta}$. Using the fact that $\pmb{x}'\in\Omega_{\epsilon,\eta}$, we  have, if $n$ is big enough, $(x'_{n-1}|x'_n)>(l-\epsilon)n$ ; and thus, using Proposition \ref{propriété de la métrique sur le bord},
$$d_a(x'_{n-1},x'_n) \leq \lambda a^{-n(l-\epsilon)}\ .$$
So for $m>n$,
$$d_a(x'_{n},x'_m) \leq \lambda a^{-(n+1)(l-\epsilon)}\sum_{k=0}^{m-1}a^{-k(l-\epsilon)}\ ;$$
which yields
$$d_a(x'_{n},x'_\infty) \leq \frac{\lambda a^{-(n+1)(l-\epsilon)}}{1-a^{-(l-\epsilon)}}\ .$$
Since $\pmb{x}$ clearly belongs to $C_{\pmb{x}}^n\cap\Omega_{\epsilon,\eta}$, we get
$$d_a(x_\infty,x'_\infty)\leq\frac{2\lambda a^{-(n+1)(l-\epsilon)}}{1-a^{-(l-\epsilon)}}\ .$$
\end{proof}
\noindent The previous lemma gives us a bound from below of the $\nu$-measure of the balls in $\partial X$ and so we obtain :
\begin{Prop}\label{majoration de la dimension ponctuelle par h/l}
For $\nu$-almost all $\xi$,
$$\overline{dim}_{P}\nu(\xi)\leq\frac{1}{\log a}\frac{h}{l}\ .$$
\end{Prop}
\begin{proof}
Let $\pmb{x}$ be an element in $\Lambda_{\epsilon,\eta}$. We have
$$\overline{dim}_{P}\nu(x_{\infty})=
        \limsup_{r\rightarrow 0}\frac{\log \nu B(x_{\infty},r)}{\log r}
               =\limsup_{r\rightarrow 0}\frac{\log \mathbb{P} D(x_{\infty},r)}{\log r}\ ;$$
and replacing $r$ by $\frac{2\lambda a^{-(n+1)(l-\epsilon)}}{1-a^{-(l-\epsilon)}}$,
$$\overline{dim}_{P}\nu(x_{\infty})
     =\limsup_{n\rightarrow \infty} \frac{\log \mathbb{P} D\Big(x_{\infty},\frac{2\lambda a^{-(l-\epsilon)(n+1)}}{1-a^{-(l-\epsilon)}}\Big)}
{-\log a\ (n+1)(l-\epsilon)}\ .$$
Using  Lemma \ref{not too far from x{infini} (Lem)}, we get
$$\overline{dim}_{P}\nu(x_{\infty})\leq  \limsup_{n\rightarrow \infty} \frac{\log \mathbb{P} \big(C_{\pmb{x}}^n\cap\Omega_{\epsilon,\eta}\big)}
{-\log a\ (n+1)(l-\epsilon)}\ ;$$
then using  Lemma \ref{Kaimanov. (Lem)}:
$$\overline{dim}_{P}\nu(x_{\infty})\leq  \limsup_{n\rightarrow \infty} \frac{-\log \mu^n(x_n)}
{\log a\ (n+1)(l-\epsilon)}= \frac{h}{\log a\ (l-\epsilon)}\ .$$
This being true for each $\epsilon>0$ on a set of measure $1-2\eta$ for each $\eta$, it proves the announced result.
\end{proof}

\section{An integral formula for the rate of escape}\label{Une formule intégrale pour la vitesse de fuite}

\subsection{Busemann functions}
We recall the definition of Busemann functions, which we will use to estimate the rate of escape. Let $\sigma$ be a geodesic ray. For each $x$ in $X$, the function $t\mapsto d(x,\sigma(t))-t$ is decreasing and bounded (it's just a consequence of the triangle inequality). So we can define the Busemann function associated to a geodesic ray $\sigma$ as :
$$f_{\sigma}(x)=\lim_{t\rightarrow\infty}d(x,{\sigma}(t))-t \ .$$
Then we define a cocycle on $X^2$ by :
$$\beta_{\sigma}(x,y)=f_{\sigma}(y)-f_{\sigma}(x) \ .$$
\begin{Rem}\label{cas où les rayons équivalents sont asymptotiques}
Assume in addition that the space $X$ has the following property $(P)$ : for two equivalent rays ${\sigma}_1$ and ${\sigma}_2$, there exists $T_0$ such that $\lim_{t\rightarrow\infty}d({\sigma}_1(t),{\sigma}_2(t+T_0))=0 \ .$  We then have, if ${\sigma}_1 \sim {\sigma}_2$,
\begin{equation}\label{égalité des cocycle (eq)}
\beta_{{\sigma}_1}(x,y)=\beta_{{\sigma}_2}(x,y) \ ;
\end{equation}
and so the function $\beta_{\sigma} $ depends only on the endpoint of the ray ${\sigma}$. This allows us to define a cocycle on $\partial X$, called the  Busemann cocycle, by 
$$\beta_\xi(x,y)=\beta_{\sigma}(x,y) \ ;$$
${\sigma}$ being a ray with endpoint $\xi$.\\
The hyperbolic halfplane $\mathbb{H}^2$, trees, an d more generally CAT(-1) spaces have this property.% ; but we shall not assume it is satisfied in the sequel.
\end{Rem}
\noindent In the general case, this  construction fails because  if ${\sigma}_1$ and ${\sigma}_2$ have the same endpoint $\xi$, which means that there exists a constant $D$ such that for all $t$,
$d(\sigma_1(t),\sigma_2(t))\leq D$, then the equality (\ref{égalité des cocycle (eq)}) is not satisfied but we only have
$$\big|\beta_{{\sigma}_1}(x,y)-\beta_{{\sigma}_2}(x,y)\big|\leq 2D \ .$$
However, two geodesic rays with the same 
point at infinity have the following property : there exists $T_0$ such that
$\lim_{t\rightarrow\infty}d({\sigma}_1(t),{\sigma}_2(t+T_0))\leq 16 \delta$ (\cite{GH}). This implies 
\begin{Lem}\label{lemCoo1}
Let ${\sigma}_1$ and ${\sigma}_2$ be two geodesic rays with the same 
point at infinity $\xi$. Then for any $x$ and $y$ in $X$, we have
$$\big|\beta_{{\sigma}_1}(x,y)-\beta_{{\sigma}_2}(x,y)\big|\leq C_1 \ ;$$
where $C_1$ is a constant relied only on $\delta$.
\end{Lem}
\noindent This allows us to adopt the following definition : for $\xi$ in $\partial X$ and $x$, $y$ in $X$,
$$\beta_{\xi}(x,y)=\sup\{\beta_{\sigma}(x,y)\} \ ;$$
where the supremum is taken on all rays with endpoint $\xi$. In particular, if ${\sigma}$ is a geodesic ray such that ${\sigma}(\infty)=\xi$, then
\begin{equation}\label{def.coc.}
\big|\beta_{\xi}(x,y)-\beta_{\sigma}(x,y)\big|\leq C_1 \ .
\end{equation}
 This will not be a  cocycle but we have
$$\big|\beta_{\xi}(x,y)-\big(\beta_{\xi}(x,z)+\beta_{\xi}(z,y)\big)\big|\leq 3C_1 \ .$$

\noindent We will need the following lemma (\cite{Coo2}) :
\begin{Lem}\label{lemCoo2}
Let $\xi$ be a point in $\partial X$, ${\sigma}$ be a ray with endpoint $\xi$, and $x_1$, $x_2$ be two points in $X$. Then there exists a neighborhood $V$
of $\xi$ such that if $y\in X\cap V$,
$$\big|\beta_{\sigma}(x_1,x_2)-\big(d(x_2,y)-d(x_1,y)\big)\big|\leq C'\ ;$$
where  $C'$ is a constant which depends  only on $\delta$. In particular we have
$$\big|\beta_{\xi}(x_1,x_2)-\big(d(x_2,y)-d(x_1,y)\big)\big|\leq C'+C_1=C_2\ .$$
\end{Lem}
\begin{Rem}\label{}
In order to estimate the rate of escape, we introduce some constants $C_i$, notation we shall keep in what follows.
\end{Rem}

\subsection{Rate of escape}
We shall get a formula for the rate of escape of a random walk in term of  this Busemann "quasi-cocycle". Let $\mu$ be a measure supported by a subgroup $G$ of the  group of isometries of $(X,d)$, with a finite first moment.
\begin{Prop}\label{formule pour la vitesse de fuite}
We assume that the random walk $(x_n o)$ associated to $\mu$
converges almost surely to an element $x_{\infty}$ in $\partial X $ (which we saw is the case under hypothesis adopted in part \ref{Poisson boundary and harmonic measure}); denote by $\nu$ the distribution of $x_{\infty}$. We then have
$$\Big|l(G,\mu)-\sum_g \mu(g)\int_{\partial X}\beta_{\xi}(o,g^{-1}o)d\nu(\xi)\Big|
\leq C_2\ ;$$ where $C_2$ is the constant introduced in Lemma \ref{lemCoo2}.
\end{Prop}
\begin{proof}
Writing
$$L_n=\int d(o,x_n o)d\mathbb{P}\ ;$$
we have by definition $$l(G,\mu)=\lim_n\frac{L_n}{n}\ .$$
So we have
\begin{equation}\label{liminf}
\liminf_{n}(L_{n+1}-L_n)\leq l(G,\mu)\leq \limsup_{n}(L_{n+1}-L_n)\ .
\end{equation}
Besides we have
\begin{align*}
L_{n+1}=&\int_{G} d(o,\gamma o)d\mu^{n+1}(\gamma)\\
   =&\sum_g \mu(g)\int d(o,g\gamma o)d\mu^n(\gamma)\\
   =&\sum_g \mu(g)\int d(o,g x_n o) d\mathbb{P} \ ;
\end{align*}
and so
\begin{align*}
L_{n+1}-L_n=&\sum_g \mu(g)\int\big(d(o,g x_n o)-d(o,x_n o)\big)d\mathbb{P}\\
            =&\sum_g \mu(g)\int\big(d(g^{-1}o,x_n o)-d(o,x_n o)\big)d\mathbb{P} \ .
\end{align*}
The quantity $\big(d(g^{-1}o,x_n o)-d(o,x_n o)\big)$ is bounded for every $n$ by $d(o,go)$,
which is an integrable function w.r.t. $\mu\otimes\mathbb{P}$ ; hence we can applly the Lebesgue convergence theorem, which implies :
$$\liminf_{n}(L_{n+1}-L_n)
     \geq \sum_g \mu(g)\int\liminf_{n}\big(d(g^{-1}o,x_n o)-d(o,x_n o)\big)d\mathbb{P} \ .$$
Moreover Lemma \ref{lemCoo2} gives
$$\big|\liminf_n \big(d(g^{-1}o,x_n o)-d(o,x_n o)\big)-\beta_{x_{\infty}}(o,g^{-1}o)\big|
                                                     \leq C_2\ ;$$
from which we deduce
$$l(G,\mu)\geq \sum_g \mu(g)\int \beta_{x_{\infty}}(o,g^{-1}o)d\mathbb{P} - C_2\ .$$
We do the same for the upper limit.
\end{proof}
\begin{Rem}\label{}
If $(X,d)$ has property $(P)$ (see Remark \ref{cas où les rayons équivalents sont asymptotiques}) the previous formula is exact.
\end{Rem}

\section{Construction}\label{Construction}
In this section $(X,d)$ is a $\delta$-hyperbolic space, $G$ a non-elementary subgroup of the group of isometries on $X$ acting properly discontinuously. 
Our goal is to construct a random walk such that the associated $h/l$ be arbitrarily small. For this purpose, we are going to construct a sequence $(\mu_k)$ of probability measures such that 
  $$\lim_{k\rightarrow \infty}\frac{h(G,\mu_k)}{l(G,\mu_k)}=0\ .$$

\subsection{Introduction of the sequence $(\mu_k)$}\label{Introduction de la suite}
 We fix  a measure $\mu$ on $G$ with a finite first moment, and such that its support generates $G$ as a semigroup.
We fix also a hyperbolic element (see the definition \ref{élément hyperbolique}) $\gamma_0$ in $G$ ; we denote $\gamma_k = \gamma_0^k o$ and $\gamma_{\pm} = \lim_{k\rightarrow \pm\infty}\gamma_k$. Then for each $k\geq 0$ we take
$$ \mu_k = \frac{1}{2}\mu +
                     \frac{1}{4}\big(\delta_{\gamma_{0}^{k}}+\delta_{\gamma_{0}^{-k}}\big) \ .$$
Remark that each $\mu_k$ satisfies the hypothesis of the part \ref{Random walk}. We denote by $\nu_k$ the harmonic measure  associated with each $\mu_k$. 

\begin{Prop}\label{entropies bornées}
The quantity $h(G,\mu_k)$ is bounded by a constant  which is not relied on $k$.
\end{Prop}
\begin{proof}
We know that $h(G,\mu_k)$ bounded from above by $H(\mu_k)$ and we have
\begin{align*}
H(\mu_k)&=-\frac{1}{2}\sum_g\mu(g)\log\Big(\frac{1}{2}\mu(g)
                                           +\frac{1}{4}\big(\delta_{\gamma_k}(g)
                                           +\delta_{\gamma_{-k}}(g)\big)\Big)\\
             &\qquad \qquad\qquad\qquad\qquad               -\frac{1}{4}\log\Big(\frac{1}{2}\mu(\gamma_k)+\frac{1}{4}\Big)
                                                 -\frac{1}{4}\log\Big(\frac{1}{2}\mu(\gamma_{-k})+\frac{1}{4}\Big)\\
           &\leq -\frac{1}{2}\sum_g \mu(g)\log\Big(\frac{1}{2}\mu(g)\Big)-\frac{1}{2}\log\Big(\frac{1}{4}\Big)\leq\frac{3\log 2}{2}+\frac{1}{2}H(\mu)\ . 
\end{align*}
\end{proof}
\noindent Hence it just remains  to show that $l(G,\mu_k)$ goes to infinity. By using the formula of   Proposition \ref{formule pour la vitesse de fuite}  we get :
$$l(G,\mu_k)\geq \sum_g \mu_k(g)\int_{\partial X}\beta_{\xi}(o,g^{-1}o)d\nu_k(\xi)-C_2\ .$$
But
\begin{align*}\sum_g \mu_k(g)\int_{\partial X}\beta_{\xi}(o,g^{-1}o)d\nu_k(\xi) =&
  \frac{1}{2}\sum_g \mu(g)\int_{\partial X}\beta_{\xi}(o,g^{-1}o)d\nu_k(\xi)  \\
  & \qquad +\frac{1}{4}\int_{\partial X}
              \big[\beta_{\xi}(o,\gamma_{k})+\beta_{\xi}(o,\gamma_{-k})\big]d\nu_k(\xi)\ ;
\end{align*}
and since $\big|\beta_{\xi}(o,g^{-1}o)\big|\leq d(o,go)$, the first element in this sum is bounded in absolute value by $L(\mu)/2$ ; then, writing $C_3=C_2+L(\mu)/2$, we get
\begin{equation}\label{première estimation de la vitesse}
l(G,\mu_k)\geq
\frac{1}{4}\int_{\partial X}
       \big[\beta_{\xi}(o,\gamma_{k})+\beta_{\xi}(o,\gamma_{-k})\big]d\nu_k(\xi)-C_3\ .
\end{equation}
We are now going to estimate the quantity in the square brackets in the previous equation.

\subsection{Estimate of $\big[\beta_{\xi}(o,\gamma_{k})+\beta_{\xi}(o,\gamma_{-k})\big]$}\label{Estimation de []}

We first show a convexity inequality based on the following property of hyperbolic spaces. 
\begin{Prop}\label{}(\cite{CDP})
In the $\delta$-hyperbolic space $(X,d)$, The metric has the following quasi-convexity property : let $x_1$ and $x_2$ be two points in $X$ and $s:[0,1]\rightarrow X$ a constant speed parametrization of a segment joining $x_1$ and $x_2$. If $y$ is an other point in $X$, then we have, for all $t$ in $[0,1]$,
$$d(y, s(t))\leq t d(y,x_1)+(1-t)d(y,x_2)+4\delta \ .$$
\end{Prop}

\noindent We deduce the following property :
\begin{Prop}\label{inégalité de convexité}
There exists a constant $C_4$ depending only on $\delta$ and $\gamma_0$ such that for every $k\geq 0$ and $\xi\in\partial X$,
$$\beta_{\xi}(o,\gamma_{k})+\beta_{\xi}(o,\gamma_{-k})\geq -C_4 .$$
\end{Prop}
\begin{proof}
Let $\sigma$ be a geodesic joining $\gamma_{+}$ and $\gamma_{-}$, $K$ the Hausdorff distance between $\sigma$ and the quasi geodesic associated with $\gamma_0$ (see Proposition \ref{quasi geodesique}). Choose $o'$ on $\sigma$ such that $d(o,o')\leq K+1$, $\gamma_{\pm k}'$ on $\sigma$ such that 
$d(\gamma_{\pm k},\gamma_{\pm k}')\leq K+1$. Denote by $m_k'$ the middle point of the segment $[\gamma_{ k}',\gamma_{- k}']$.
\begin{figure}[htbp]
  \centering
\includegraphics[width=8cm]{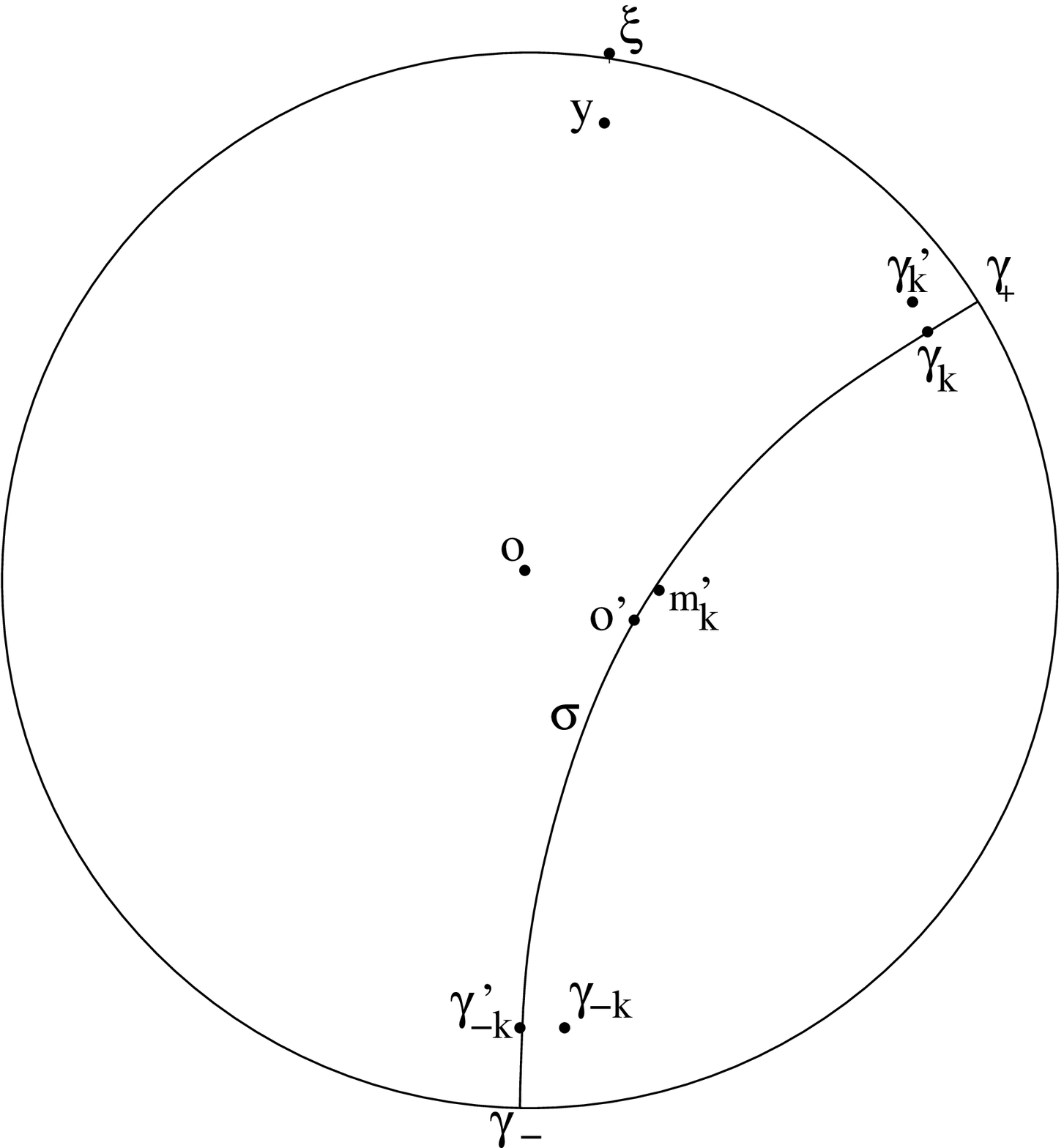}  
  \label{dessin}
\end{figure} 
For each $y$ in $X$ we have 
$$d(y,\gamma_{k})-d(y,o)+d(y,\gamma_{-k})-d(y,o)\geq d(y,\gamma_{k}')-d(y,o')+d(y,\gamma_{-k}')-d(y,o')-4(K+1)\ ;$$
and
$$d(y,\gamma_{k}')-d(y,o')+d(y,\gamma_{-k}')-d(y,o')\geq d(y,\gamma_{k}')+d(y,\gamma_{-k}')-2d(y,m_k')-2d(o',m_k')\ .$$
But using the convexity property of the distance we get
$$d(y,\gamma_{k}')+d(y,\gamma_{-k}')-2d(y,m_k')\geq -8\delta\ ;$$
and since  $d(o,\gamma_{k})=d(o,\gamma_{-k})$, 
\begin{align*}
\big|d(o',\gamma_{k}')-d(o',\gamma_{- k}')\big| &=
\Big|\big[d(o',\gamma_{k}')-d(o,\gamma_{k}')+d(o,\gamma_{k}')-d(o,\gamma_{k})\big]\\
&\qquad
-\big[d(o',\gamma_{-k}')-d(o,\gamma_{-k}')+d(o,\gamma_{-k}')-d(o,\gamma_{-k})\big]\Big|\\
 &\leq 2 d(o,o')+ d(\gamma_{ k}',\gamma_{ k})+ d(\gamma_{ -k}',\gamma_{- k})\leq 4(K+1)\ ;
\end{align*}  
which yields 
$$ d(o',m_k')\leq 2(K+1)\ .$$
So we get $$d(y,\gamma_{k})-d(y,o)+d(y,\gamma_{-k})-d(y,o)\geq -4(K+1) -8\delta -4(K+1)\ ; $$
 and using Lemma \ref{lemCoo2} allows us to conclude, denoting $8(K+1)+8\delta+C_2$ by $C_4$.
\end{proof}
\noindent We are now going to show that on a "big enough" set (w.r.t. the measures $\nu_k$) of the boundary,  $\big[\beta_{\xi}(o,\gamma_{k})+\beta_{\xi}(o,\gamma_{-k})\big]$ goes to $+\infty$ when $k$ does.\\
We have,  if $y$ is in a neigbourhood of $\xi$, again using Lemma \ref{lemCoo2},
$$\beta_{\xi}(o,\gamma_{k})+\beta_{\xi}(o,\gamma_{-k})
              \geq d(y,\gamma_{k})-d(y,o)+d(y,\gamma_{-k})-d(y,o)-C_2 \ ;$$
hence
\begin{equation}\label{formule pour []}
\beta_{\xi}(o,\gamma_{k})+\beta_{\xi}(o,\gamma_{-k})
              \geq d(o,\gamma_{k})+d(o,\gamma_{-k})-2(y|\gamma_{k})-2(y|\gamma_{-k})-C_2 \ ;
\end{equation}
so if the two Gromov products in this sum are bounded(which is the case if $\xi \notin \{\gamma_+,\gamma_-\}$), we then have
$$\lim_{k\rightarrow\infty}\beta_{\xi}(o,\gamma_{k})+\beta_{\xi}(o,\gamma_{-k})=+\infty \ .$$

\noindent However, in order to conclude using formula (\ref{première estimation de la vitesse}), we need that the quantity $\big[\beta_{\xi}(o,\gamma_{k})+\beta_{\xi}(o,\gamma_{-k})\big]$ goes to infinity uniformly on a set whose $\nu_k$-measure remains greater than a strictly positive constant. 
\begin{Lem}\label{majoration du produit de Gromov}
Let $U$ be a neigbourhood  of $\gamma_+$ in $\partial X$. Then there  exists a constant $C>0$ and an integer $K$ such that for every $k\geq K$ and every $\xi\notin U$
$$(\xi|\gamma_k)\leq C \ .$$
\end{Lem}
\begin{proof}
We saw that the topology on $\partial X$ was defined by neighborhoods of the type $\{\xi' : (\xi|\xi')\geq D'\}$. Let $C'$ be a constant such that if $(\xi|\gamma_+)\geq C'$, then $\xi\in U$. Take a $K$  such that if $k\geq K$, then $(\gamma_k|\gamma_+)\geq C' + \delta$. Let then be $\xi\notin U$ ; if we had $(\gamma_k|\xi)\geq C' + \delta$, we would have
$$(\xi|\gamma_+)\geq \min \big\{(\gamma_k|\gamma_+),(\gamma_k|\xi)\big\}-\delta \geq C' \ ;$$
so we deduce the result taking $C=C'+\delta$.
\end{proof}

\noindent We saw in Proposition \ref{la mesure harm charge les ouverts} how to show that the harmonic measure of an open set  was strictly positive by using the minimality of the limit set through Proposition \ref{conséquence de la minimalité de l'ensemble limite}. We are now going to show in the same way that we can bound from below uniformly w.r.t. $k$ the $\nu_k$-measure of an open set.
\begin{Prop}\label{nu_k mesure des ouverts}
Let $U$ be an open set meeting $\Lambda$. There exist a constant  $\alpha>0$ (relied on $U$ and $\mu$) such that for every $k$,
$$\nu_k(U)\geq \alpha \ .$$
\end{Prop}
\begin{proof}
We take $\gamma_1,\cdots ,\gamma_r \in G$ like in Proposition \ref{conséquence de la minimalité de l'ensemble limite} :
$$\Lambda \subset \bigcup_{1\leq i \leq r}\gamma_i U \ ;$$
and $s$ such that $\gamma_1,\cdots ,\gamma_r \in supp(\mu^s)$.
Now we remark the following : we have a sort of uniform stationnarity of the $\nu_k$. Indeed as $\nu_k=\mu_k * \nu_k$ and for every $k$ $\mu_k\geq 1/2 \mu$,we have
$$\nu_k \geq \frac{1}{2^s}\mu^s * \nu_k \ ;$$
and so
\begin{align*}
\nu_k(U)&\geq \frac{1}{2^s}\sum_{g\in supp(\mu^s)}\mu^s(g)\nu_k(gU)\\
        &\geq \frac{1}{2^s}\sum_{i}\mu^s(\gamma_i)\nu_k(\gamma_i U)\\
        &\geq \frac{1}{2^s} \min_{i}\mu^s(\gamma_i)\ .
\end{align*}
\end{proof}
\noindent Since $G$ is assumed to be non-elementary, $\Lambda$ contains an element wich is distinct from $\gamma_+$ and $\gamma_-$ ; so we can fix an open set $U$ meeting $\Lambda$ and which  doesn't contains $\gamma_+$ neither $\gamma_-$. Using Lemma  \ref{majoration du produit de Gromov} with $\partial X\setminus U$ as a neighborhood of $\gamma_+$ (and $\gamma_-$), we take $C_5$ and $K$ such that if $\xi\in U$ and $k\geq K$,
$$(\xi|\gamma_{\pm k})\leq C_5 \ .$$
 Using formula (\ref{formule pour []}), we get, for $k\geq K$ and $\xi\in U$ :
$$\beta_{\xi}(o,\gamma_{k})+\beta_{\xi}(o,\gamma_{-k})\geq
     2d(o,\gamma_{k})-4 C_5-C_2\ .$$
Besides formula (\ref{première estimation de la vitesse}) and  Proposition \ref{inégalité de convexité} give us 
$$l(G,\mu_k)\geq
\frac{1}{4}\int_{U}
       \big[\beta_{\xi}(o,\gamma_{k})+\beta_{\xi}(o,\gamma_{-k})\big]d\nu_k(\xi)
                                   -\frac{C_4}{4}+C_3\ ;$$
and now we use the previous proposition to get 
$$l(G,\mu_k)\geq \frac{\alpha}{4}\big[2d(o,\gamma_{k})-4 C_5-C_2\big]-\frac{C_4}{4}-C_3\ .$$
Since the first element of the right member in previous equation goes to infinity, this gives us :
\begin{Prop}\label{convergence des vitesses de fuite vers l'infini}
With the notations previously adopted,
$$\lim_{k\rightarrow\infty}\frac{h(G,\mu_k)}{l(G,\mu_k)}=0 \ .$$
\end{Prop}

\section{Main result}\label{Résultat principal}
Putting together  Proposition \ref{convergence des vitesses de fuite vers l'infini} of the previous section and  Proposition \ref{majoration de la dimension ponctuelle par h/l} of section \ref{Lien entre les quantité asymptotiques et la dimension de Hausdorff de la mesure harmonique}, we get the following result :

\begin{theo}\label{gros théorème}
Let $G$ be a subgroup of the group of isometries of a hyperbolic space $(X,d)$ acting properly discontinuously and which is not elementary. For every $\epsilon>0$ there exists on $G$ a symmetric probability measure $\mu$ with a finite first moment, whose support generates $G$, and which has the following property: the pointwise dimension of the harmonic measure $\nu$ associated with $(G,\mu)$ is $\nu$-almost surely smaller than $\epsilon$.
\end{theo}
%\begin{Rem}\label{degenerated}In what extense can we say the measures we have constructed are not "degenerated" ? We demand of course that the support of $\mu$ generates $G$. But it is not sufficient : in the case of the free group with two generators $a$ and $b$, if we take $\mu(a^{-1})=\mu(b^{-1})=\mu(b)=\epsilon$ and $\mu(a)=1-3\epsilon$, the dimension of the harmonic measure become small when $\epsilon$ goes to zero, but to some extent this exemple can be considered as "degenerated". We avoid this problem in our construction because the support of $\mu$ can be taken symmetric~;~ remark that for the free group, when $\mu$ is taken symmetric on the generators, we can't have an arbitrarily small dimension for the harmonic measure. Moreover, a measure $\mu_0$ (with associated harmonic measure having strictly positive dimension) being fixed, we can assume that our measure $\mu$ is greater than $\mu_0/2$.
%\end{Rem}

\noindent Now we restrict ourself to the case where $(X,d)$ is the Cayley graph of $G$. In this case, in view of Propositions \ref{majoration de la dim de H par la dim ponctuelle} and \ref{multiplicité finie du bord d'un groupe hyperbolique}, Theorem \ref{gros théorème} implies that the Hausdorff dimension of the harmonic measure can be strictly lower than the dimension  of the boundary, whose value is $e_a(G)$. And in view of Remark \ref{la dimension caractérise le type}, it implies that the harmonic and Hausdorff measures are not equivalent; since both are ergodic they are singular. So we get the following result on the type of the harmonic measure~:
\begin{Cor}\label{corollaire sur les groupes (Cor)}
In the case where $X$ is the Cayley graph of a hyperbolic group $G$, there exist on $G$ a symmetric probability measure $\mu$ with a finite first moment, whose support generates $G$, and such that the associated harmonic measure and the Hausdorff measure on $\partial G$ are mutually singular.
\end{Cor}\label{}
\vspace{0.2cm}
\noindent Let us finally note that concerning the asymptotic behaviour of the sequence $(\nu_k)$, we can prove that for every subsequence of $(\nu_k)$ which converges weakly to a certain measure $\nu$, all points of the orbit of $\gamma_+$ and $\gamma_-$ under $G$ are atoms for $\nu$. In the case of a finitely generated free group, the sequence $(\nu_k)$ even converges weakly to an atomic measure suported by the orbit of $\gamma_+$ and $\gamma_-$ under $G$.

\vspace{1cm}

\noindent \begin{bf}Acknowledgments.\end{bf} I would like to thank my supervisor Vadim Kaimanovich for his advices during the redaction of this paper.

%\newpage

\end{document}